\documentclass{article} 
\usepackage{iclr2027_conference,times}

\usepackage{amsmath,amsfonts,bm}

\def\1{\bm{1}}

\DeclareMathAlphabet{\mathsfit}{\encodingdefault}{\sfdefault}{m}{sl}
\SetMathAlphabet{\mathsfit}{bold}{\encodingdefault}{\sfdefault}{bx}{n}

\usepackage{hyperref}
\usepackage{url}
\usepackage{amsmath,amssymb}   
\usepackage{xcolor}            
\usepackage{graphicx}
\usepackage{subcaption}
\usepackage{float}
\usepackage[title]{appendix}

\newcommand{\uhat}{\hat{\mathbf{u}}}
\newcommand{\arf}{a_{\mathrm{RF}}}

\title{Amortized Generative modeling of Invariant Measures Across a Family of Kuramoto-Sivashinsky PDEs}

\author{Arnab Roy \& Tobin A. Driscoll\\ 
Department of Mathematical Sciences\\
University of Delaware\\
Newark, DE 19716, USA \\
\texttt{\{arnabroy, driscoll\}@udel.edu} \\
}

\iclrfinalcopy 
\begin{document}

\maketitle

\begin{abstract}
The long-run statistics of chaotic dissipative partial differential equations are described by invariant measures, and when the equation carries a parameter these measures form a family whose members can differ in kind, from steady states to sustained spatiotemporal chaos. We ask whether a single generative model can cover such a family. We train a single conditional diffusion model once, on states of the Kuramoto–Sivashinsky equation at 160 values of its hyperviscosity, and evaluate it against held-out data at interleaved parameter values not seen in training, using a classifier two-sample test supported by spectral, geometric and tail statistics. Over most of the family, agreement improves as the dynamics become more chaotic.
We identify a structural cause: the weakly chaotic members have invariant measures concentrated on low-dimensional sets, which a model with full-support output cannot represent, and the target approaches full dimensionality as chaos increases. Restricting training to the strongly chaotic members reduces the classifier's excess over chance on those members by roughly 40\%, and this restricted conditional model matches or exceeds models trained at single parameter values. On those members, generated samples reproduce energy spectra to within a resampling floor, with a small systematic deficit at the largest scales, and match the energy balance of held-out data.
\end{abstract}

\section{Introduction}
\label{sec:intro}

Chaotic dissipative partial differential equations possess invariant measures
that describe their long-run statistics
\citep{foias2001}
For many purposes it is those statistics, rather than any individual
trajectory, that are the object of interest: in weather and climate modeling
the distribution a system settles into and the skill of a particular forecast
are distinct goals, pursued with different tools
\citep{price2025, meuer2026}. Drawing samples from such a measure is a
generative modeling problem, and diffusion models have recently been shown
to solve it for a chaotic ordinary differential equation
\citep{finn2024}, alongside broader use of generative diffusion for
forecasting and ensemble generation in physical systems \citep{nai2025}.

Physical systems, however, come as parametrized families rather than as
isolated instances, and the parameter can change the dynamics qualitatively.
The Kuramoto--Sivashinsky equation is the standard illustration: as its
domain is enlarged the dynamics pass from a stable fixed point through
traveling waves and periodic orbits to sustained spatiotemporal chaos
\citep{papageorgiou1991}, so that a single scalar sweeps out a family of
systems whose invariant measures differ in kind and not merely in detail.
Fitting each member separately is the obvious approach, and its cost grows
with the number of parameter values one cares about. Parameter-conditioned
networks have been used to emulate dynamics across such families
\citep{shokar2025}, and we ask whether a single model can reproduce the
invariant \emph{measures} of a whole family.

We study the Kuramoto--Sivashinsky equation on a fixed periodic domain with
hyperviscosity $\nu$ as the parameter, which is equivalent to the more familiar
formulation in which the domain length varies. We train one conditional
diffusion model on $\ln \nu$, once, on data pooled across $160$ values of the
parameter spanning steady through strongly chaotic dynamics, and ask how well
$p_\theta(\cdot \mid \ln\nu)$ reproduces the invariant measure $\mu_\nu$ at
parameter values held out from training. Accuracy is measured by a classifier
two-sample test \citep{lopez-paz2017}, supported by spectral, geometric and
tail statistics, and by a check that the generated samples satisfy the energy
balance the dynamics impose.

Our findings are the following.
\begin{itemize}\setlength{\itemsep}{2pt}
  \item \textbf{Amortization is free, and sometimes better than free.} A
  single conditional model matches or exceeds the accuracy of models trained
  at one parameter value each; at a moderately chaotic member the conditional
  model, trained across the whole family, fits that member more closely than
  a model trained on it alone.
  \item \textbf{Accuracy improves as the dynamics become more chaotic, for a
  structural reason.} The classifier accuracy falls from complete
  separability toward the chance floor as the leading Lyapunov exponent
  grows. The weakly chaotic members are not fit badly because chaos is easy
  and regularity is hard, but because their invariant measures concentrate on
  low-dimensional sets that a model with full-support output cannot
  represent; the target becomes representable as the dynamics become chaotic.
  \item \textbf{A family helps only when its members are of the same kind.}
  Restricting training to the strongly chaotic sub-family reduces the excess
  of classifier accuracy over chance by roughly $40\%$ on those same members,
  so capacity spent on targets the model cannot represent competes with
  targets it can.
  \item \textbf{The generated measures respect the constraints the dynamics
  impose.} Energy spectra agree with held-out data to within a resampling
  floor over most of the range, and the effective hyperviscosity implied by the
  generated samples matches that of the held-out samples to a median of
  $0.2\%$. The one systematic discrepancy is a deficit of energy at the
  largest scales, which follows the neutral wavenumber across the family.
\end{itemize}

The parameter values used for testing interleave those used for training, so
these results characterise interpolation within the family rather than
extrapolation beyond it. We make no claim about the relative merits of
diffusion models and other generative families for this problem; the
comparisons here are within the family and against held-out data, loosely following
the evaluation protocol of \citet{finn2024}.

\section{Related work}
\label{sec:related}

The closest precedent to this work is \citet{finn2024}, who train an
unconditional diffusion model on the climatology of the Lorenz-63 system and evaluate it with Hellinger distances, nearest-neighbor statistics and tail coverage; our evaluation begins from theirs, retaining several of those tests and adding others suited to a spatially extended system.
Two things differ.
Their target is a three-dimensional ordinary differential equation, whereas ours is
a partial differential equation whose truncated state has $38$ dimensions and
whose attractor dimension exceeds $20$ at the chaotic end of our family; and
their model is fitted to a single system, whereas ours is conditioned on a
parameter and covers a family. Generative diffusion has also been applied
widely to physical systems for forecasting and ensemble generation
\citep{price2025, meuer2026, nai2025}. Those models are conditioned on a
current state and sample a subsequent one, which is a different object from
the invariant measure: they model the transition kernel of the dynamics,
while we model the distribution the dynamics settle into.

The Kuramoto--Sivashinsky equation is a standard testbed for data-driven
models of spatiotemporal chaos.
Reservoir computing \citep{pathak2018} and recurrent architectures \citep{vlachas2020} have both been demonstrated on it, the latter across a range of domain lengths.
This line of work learns the dynamics, in the form of a map from the current state to a future one, and is assessed by how long a prediction remains close to the truth, measured in Lyapunov times.
Our object is different: we learn the invariant measure and
never predict a trajectory, so forecast horizon is not a meaningful criterion
for what we do, and distributional agreement is not a meaningful criterion
for a forecaster over short horizons.

Conditioning a single network on the parameters of a partial differential
equation, so that one model covers a family rather than an instance, has been
used for emulation of chaotic dynamics \citep{shokar2025}; the construction
here is the analogue for invariant measures. The underlying idea is
amortization: a network trained across many instances replaces per-instance
fitting at deployment, as in amortized Bayesian inference, whose software we
use for our implementation \citep{radev2023}.

Finally, the mechanism we identify for the failure of our models in 
weakly chaotic members is related to the known difficulties of score-based
generative models on targets concentrated near low-dimensional sets
\citep{pidstrigach2022}. What we report is an instance of that difficulty
arising from the physics of a system rather than from the geometry of a
dataset: it is the dynamics, through the parameter, that determine whether
the target measure is representable.

\section{Problem setup}
\label{sec:setup}

\subsection{The Kuramoto--Sivashinsky equation}
\label{sec:ks}

We aim to reproduce solutions of the Kuramoto--Sivashinsky (KS) equation in the form
\begin{equation}
\label{eq:ks}
u_t + u u_x + u_{xx} + \nu u_{xxxx} = 0,
\qquad x \in [0, 2\pi), \quad t \ge 0,
\end{equation}
with periodic boundary conditions and a single parameter $\nu \ge 0$. This
equation was derived independently as a model of phase turbulence in
reaction--diffusion systems \citep{kuramoto1976} and of flame-front
instability in laminar combustion \citep{sivashinsky1977}. It is among the
simplest partial differential equations to exhibit sustained spatiotemporal
chaos, and is for that reason a standard testbed for data-driven models of
chaotic dynamics.

The hyperviscosity parameter $\nu$ controls how complicated the dynamics are. Linearizing
\eqref{eq:ks} about $u \equiv 0$, the Fourier mode $e^{ikx}$ has growth rate
$k^2 - \nu k^4$, so modes with $k < k_* = \nu^{-1/2}$ are linearly unstable
and all higher modes are damped, while the nonlinear advective term redistributes the energy from stable to unstable modes.
Decreasing $\nu$ creates more unstable modes, and the dynamics generally grows
more complex: the system passes from a stable fixed state, through traveling
waves and periodic orbits, towards sustained chaos \citep{papageorgiou1991}.
The progression is not monotonic; however, windows of steady or periodic
behavior persist within the predominantly chaotic range.
A single scalar knob therefore tunes a family of dynamical systems of increasing complexity, making KS a natural setting for the question we ask: whether one conditional generative model, trained once, can reproduce the long-time statistics of every member of the family.

The form \eqref{eq:ks} fixes the domain and varies the hyperviscosity $\nu$. The KS literature more commonly chooses $\nu=1$ and
varies the size of the domain, writing $v_\tau + v v_y + v_{yy} + v_{yyyy} = 0$ on $y \in [0, L)$. The two forms are equivalent under a rescaling of space, time and amplitude (Appendix~\ref{app:rescaling}), with
\begin{equation}
\label{eq:L}
L = 2\pi \nu^{-1/2}, \qquad \nu = \left(2\pi/L\right)^{2}.
\end{equation}
We adopt this convention for reporting our results, but all computation is carried out for~\eqref{eq:ks} on the fixed domain
$[0,2\pi)$. In particular, the retained wavenumbers are integers throughout, and not of the form $2\pi n / L$.

\subsection{State representation}
\label{sec:state}

Since $u(\cdot,t)$ is real and $2\pi$-periodic, it is determined by its Fourier
coefficients,
\begin{equation}
\label{eq:fourier}
u(x,t) = \sum_{k \in \mathbb{Z}} \hat{u}_k(t)\, e^{ikx},
\qquad \hat{u}_{-k} = \overline{\hat{u}_k},
\end{equation}
with the modes with $k \ge 0$ carrying all the information. Furthermore,~\eqref{eq:ks} conserves the spatial mean, so that $\hat{u}_0$ is constant in
time (Appendix~\ref{app:mean}). We consequently set $\hat{u}_0 \equiv 0$ for all time and exclude it from the state.

We retain the modes $k = 1, \dots, K$ with $K = 19$. This cutoff lies above
the neutral wavenumber $k_* = \nu^{-1/2}$ for every $\nu$ in our experiments, so
all linearly unstable modes are represented, and discarded modes carry a negligible share of the total
energy. This state truncation is coarser than the resolution of the numerical KS solver, which is
described in Appendix~\ref{app:solver}.

Separating each retained coefficient into its real and imaginary parts gives
the state vector
\begin{equation}
\label{eq:state}
\uhat =
\bigl(\operatorname{Re}\hat{u}_1, \dots, \operatorname{Re}\hat{u}_{K},\;
      \operatorname{Im}\hat{u}_1, \dots, \operatorname{Im}\hat{u}_{K}\bigr)
\in \mathbb{R}^{d}, \qquad d = 2K = 38 .
\end{equation}
Real and imaginary parts are stacked in blocks rather than interleaved. This
is the representation in which all models are trained and all distributional
comparisons are made.

\subsection{The invariant measure and the learning problem}
\label{sec:measure}

For fixed $\nu$, the semiflow generated by equation~\eqref{eq:ks} possesses a
global attractor of finite dimension
\citep{nicolaenko1985}
and the long-time statistics of
the dynamics are described by an invariant probability measure $\mu_\nu$
supported on it. We work throughout with the pushforward of this measure under
the projection onto the retained modes, so that $\mu_\nu$ is a measure on
$\mathbb{R}^{d}$ and $\uhat \sim \mu_\nu$ is meaningful. We assume the
dynamics to be ergodic with respect to $\mu_\nu$, so that time averages along
a single trajectory converge to averages against $\mu_\nu$.
Our protocol to sample from $\mu_\nu$ is to record snapshots at equally spaced times along one long trajectory, after a burn-in period to remove transients (Appendix~\ref{app:solver}).

Varying $\nu$ yields not a single measure but a one-parameter family
$\{\mu_\nu\}$. Our object of study is a single conditional model
\begin{equation}
\label{eq:target}
p_\theta(\uhat \mid c), \qquad c = \ln \nu,
\end{equation}
with one parameter vector $\theta$, trained once on data pooled across the
family, and intended to satisfy $p_\theta(\,\cdot \mid \ln\nu) \approx \mu_\nu$
for every $\nu$ in the family. The remainder of the paper is concerned with
how accurate that approximation is, and with how its accuracy varies across
the family.

\subsection{Quantifying chaos across the family}
\label{sec:lyapunov}

The parameter $\nu$ tunes the behavior of KS but is an imperfect
proxy for how complicated the dynamics are: within ranges of $\nu$ that are
predominantly chaotic, there are windows in which trajectories settle onto
fixed points or traveling waves
\citep{papageorgiou1991}.
We therefore use the leading
Lyapunov exponent $\lambda_1$, the asymptotic rate at which nearby
trajectories separate, to indicate how strongly chaotic a given member of
the family is.
We estimate $\lambda_1$ by the Benettin algorithm \citep{benettin1980},
accumulated over the whole retained trajectory, and report it in inverse time
units of~\eqref{eq:ks}; the exponent is defined in
Appendix~\ref{app:lyapunov}.


\section{Method}
\label{sec:method}

\subsection{Data}
\label{sec:data}

We instantiate the family at $336$ values of $\nu$, defined by
$\nu = \mu^{-2}$ with $\mu$ uniformly spaced on $[1.1, 16.1]$. Since
$\mu = k_* = \nu^{-1/2}$, this is a uniform grid in the neutral wavenumber
and equivalently in the domain length $L = 2\pi\mu$, which ranges from
$6.91$ to $101.16$ in steps of $0.281$; the corresponding hyperviscosities run
from $\nu \approx 0.826$ down to $\nu \approx 0.0039$. The number of linearly
unstable modes $\lfloor \mu \rfloor$ increases from $1$ to $16$ across the
family, which accordingly spans members whose trajectories are steady in time
through to strongly chaotic ones.

The values are assigned to training, test and validation sets by alternating
along the grid, giving $160$ training values, $160$ test values and $16$
validation values. Training and test values therefore interleave throughout,
and the test range lies strictly inside the training range at both ends. For
each $\nu$ we integrate one trajectory from an initial condition common to the
whole family, discard a transient, and retain the final $10{,}000$ states,
spanning $181.8$ time units; solver, time grid and
transient length are given in Appendix~\ref{app:solver}.

\subsection{Model}
\label{sec:model}

We use a diffusion model \citep{ho2020, song2021}. We
train both, conditional models on the full family and on the strongly chaotic
sub-family of Section~\ref{sec:lyapunov}, and unconditional models at
individual values of $\nu$. Full details of the diffusion formulation, the
network and the training runs are given in Appendix~\ref{app:training}.

\section{Evaluation}
\label{sec:evaluation}

All evaluation is two-sample: for each test value of $\nu$, $10{,}000$
generated states $q = \{\uhat^{\rm gen}_n\}$ are compared against
$10{,}000$ held-out states $p = \{\uhat^{\rm test}_m\}$ at the same $\nu$. 

Our primary metric is a classifier two-sample test
\citep{lopez-paz2017}. A random-forest classifier (whose settings are given in Appendix~\ref{app:metrics}) is trained to
distinguish generated from held-out states, and we report its accuracy
$\arf$ on data unseen in training. A model that reproduces $\mu_\nu$ exactly
yields $\arf = 0.5$: the two samples are statistically indistinguishable, and
no classifier can beat chance. Complete separability yields $\arf = 1$. To the extent that successive held-out states are correlated along the trajectory, it makes the classifier's task easier and our reported
accuracies conservative.

This is the only metric we use that is sensitive to the full joint
distribution on $\mathbb{R}^d$, including the dependence between modes, but it cannot reveal how the two distributions differ. The remaining metrics address particular aspects of that comparison.
\begin{description}\setlength{\itemsep}{0pt}
  \item[Energy spectrum.] The mean modal energies
  $E_k = \langle |\hat{u}_k|^2 \rangle$ of the generated and held-out samples are
  compared mode by mode, which localizes any mismatch in wavenumber. The
  per-mode log ratios $r_k = \log_{10} \bigl(E^{\mathrm{gen}}_k /
  E^{\mathrm{test}}_k\bigr)$ are summarized by an energy-weighted
  discrepancy
  \begin{equation}
  \label{eq:dw}
  D_w \;=\; \frac{\sum_k E^{\mathrm{test}}_k \, |r_k|}
                 {\sum_k E^{\mathrm{test}}_k},
  \end{equation}
    which counts each mode in proportion to the energy it carries. 
    While the minimum value is nominally zero for identical distributions, two finite samples typically give a nonzero $D_w$, so we compare it against a
  resampling floor: the held-out sample is split at
  random into halves, the halves are scored against each other, and the result
  is averaged over five independent splits. A model is considered spectrally
  indistinguishable from the truth when its $D_w$ sits at that floor. This
  metric sees per-mode energies only, and no dependence between modes.
  \item[Nearest-neighbor ratio \citep{finn2024}.] With
  $\|\cdot\|$ the Euclidean norm on $\mathbb{R}^d$, define
  \begin{equation}
  \label{eq:nn}
  d_{\mathrm{gen}} = \mathbb{E}_{x \sim q}
      \min_{z \in p} \|x - z\|^2,
  \qquad
  d_{\mathrm{test}} = \mathbb{E}_{z \sim p}
      \min_{x \in q} \|z - x\|^2.
  \end{equation}
  Thus, $d_{\mathrm{gen}}$ is large when the model places mass away from the
  attractor and $d_{\mathrm{test}}$ is large when it misses regions the
  held-out sample occupies. We report
  $\rho_{\mathrm{NN}} = d_{\mathrm{test}}/d_{\mathrm{gen}}$, which is our
  only metric that indicates the \emph{direction} of a failure; a
  perfect model gives $\rho_{\mathrm{NN}} = 1$. 
  \item[Tail coverage \citep{finn2024}.] Let $p^{(i)}_\beta$ denote the
  $\beta$-quantile of the held-out sample in coordinate $i$. For
  $\alpha = 0.01$,
  \begin{equation}
  \label{eq:pot}
  \mathrm{POT}_i = \mathbb{E}_{x \sim q}
     \Bigl[ \chi\bigl(x_i \leq p^{(i)}_{\alpha}\bigr)
          + \chi\bigl(x_i \geq p^{(i)}_{1-\alpha}\bigr) \Bigr],
  \end{equation}
  averaged over eight fixed coordinates of $\uhat$, where $\chi$ is the indicator function. A perfect model gives
  $\mathrm{POT}_i=2\alpha = 0.02$; smaller values indicate truncated tails, larger values
  inflated ones. This is our metric to watch for rare excursions.

  \item[Energy balance.] Multiplying equation~\eqref{eq:ks} by $u$ and
  integrating over the domain zeroes the advective
  term, and stationarity of $\mu_\nu$ forces its mean to vanish
  (Appendix~\ref{app:nueff}), so that
  $\sum_k k^2 E_k = \nu \sum_k k^4 E_k$. This motivates the effective
  hyperviscosity
  \begin{equation}
  \label{eq:nueff}
  \nu_{\mathrm{eff}} = \frac{\sum_k k^2 E_k}{\sum_k k^4 E_k},
  \end{equation}
  which tests whether a sample satisfies a balance the dynamics impose. We
  evaluate it on the generated and held-out samples separately and compare
  the two: a model that reproduces $\mu_\nu$ should return the same
  $\nu_{\mathrm{eff}}$ as the held-out data. Because the sum is truncated at $k = K$, that value need not equal $\nu$
  itself.
\end{description}

We note that the spectrum, tail and balance checks compare
low-dimensional functionals whose estimation does not degrade with dimension,
so their agreement with $\arf$ cannot be an artefact of classifier power
varying across the family.

\section{Results}
\label{sec:results}

\begin{figure}[t]
  \centering
  \includegraphics[width=\textwidth]{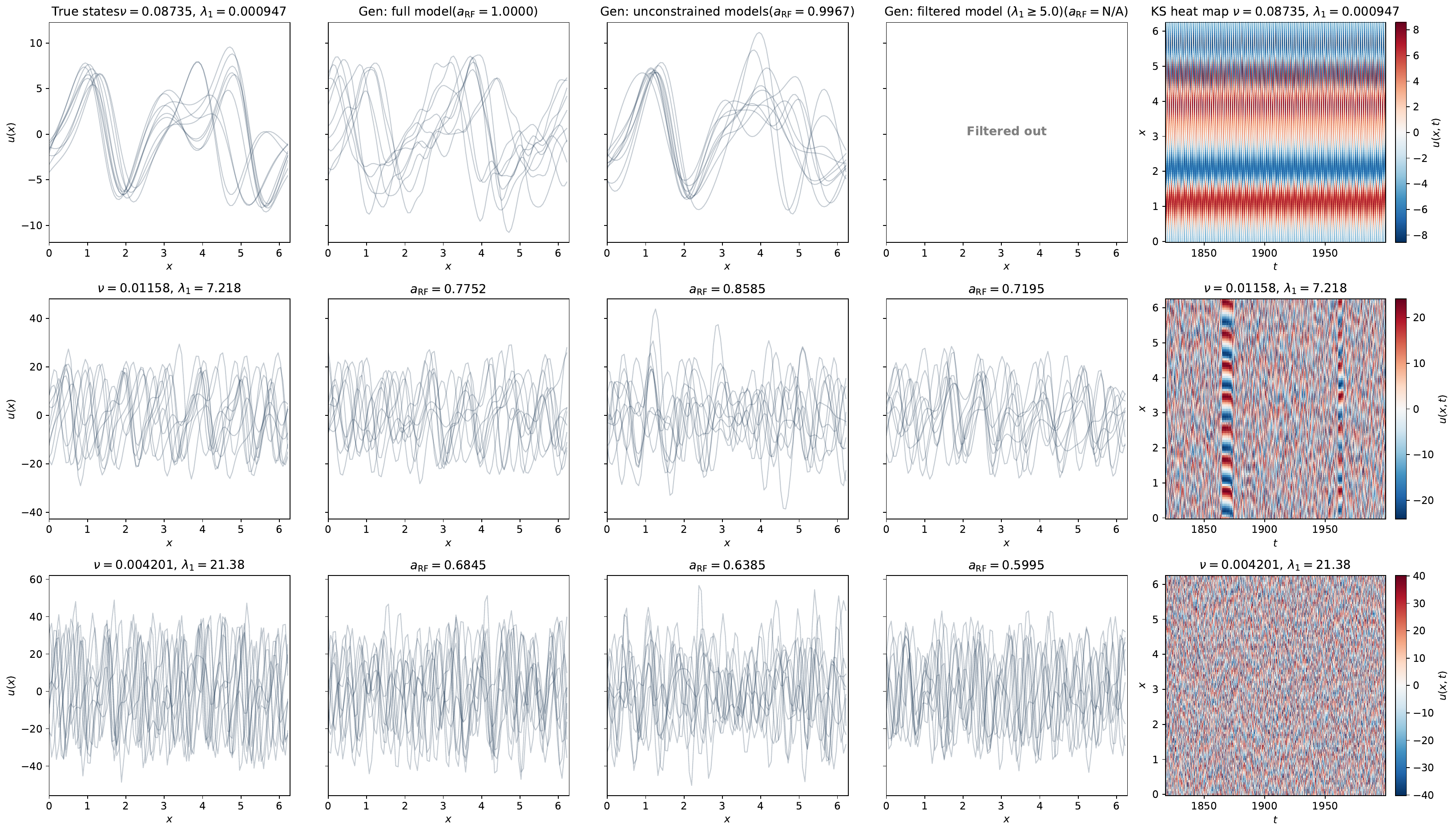}
  \caption{Generated and held-out states at three test values of $\nu$, one
per row, with $\nu$ and $\lambda_1$ given in the panel titles. Columns:
held-out states; the conditional model trained on the full family; an
unconditional model trained only at the row's $\nu$; the
conditional model trained on the $\lambda_1 \geq 5$ sub-family; the true space--time
evolution $u(x,t)$. Each panel overlays [10] independent draws reconstructed
to physical space, and panel titles give the classifier accuracy $\arf$
against the held-out sample. The weakly chaotic member (top row) lies
outside the sub-family, so the filtered model's panel is empty.}
  \label{fig:qualitative}
\end{figure}

Figure~\ref{fig:qualitative} shows examples of what the models produce. At three test
values of $\nu$, generated states are overlaid on
held-out states, with the true space--time evolution alongside. At the
moderately and strongly chaotic members the generated and held-out states are
visually indistinguishable, while the weakly chaotic member already hints at
the pattern this section quantifies. 

\subsection{Accuracy across the family}
\label{sec:headline}

\begin{figure}[t]
  \centering
  \includegraphics[width=\textwidth]{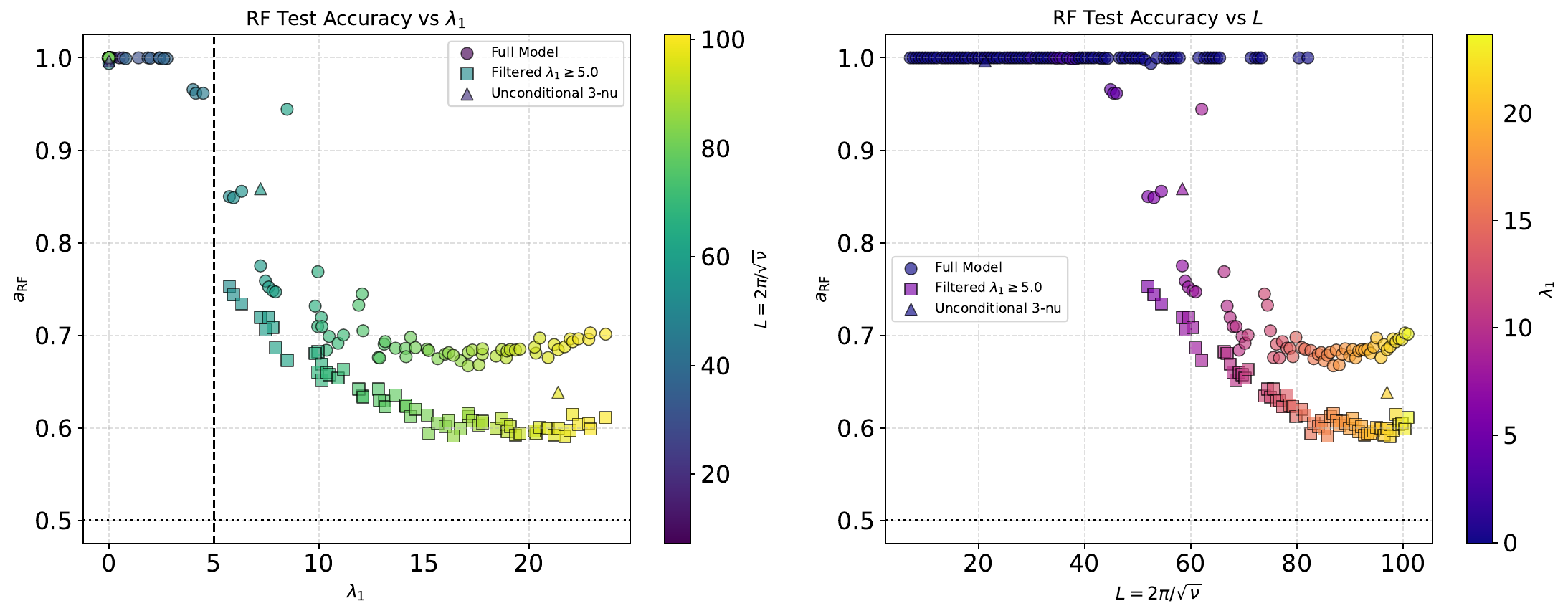}
  \caption{Classifier accuracy $\arf$ plotted against $\lambda_1$ (left) and against $L$ (right), with
  each panel colored by the other quantity. Circles: the conditional model
  trained on the full family. Squares: the conditional model trained on the
  $\lambda_1 \geq 5$ sub-family, shown at the test values it retains.
  Triangles: unconditional models trained at a single $\nu$
  (Section~\ref{sec:amortization}). The dotted line at $\arf = 0.5$ is the
  chance floor, attained when generated and held-out samples are
  indistinguishable; the dashed vertical line marks the $\lambda_1 = 5$
  threshold.}
  \label{fig:rfaccuracy}
\end{figure}

Figure~\ref{fig:rfaccuracy} reports $\arf$ for every test value of $\nu$. For
the full-family model, accuracy is essentially $1$ (i.e., worst generative performance) across the weakly chaotic
members and falls toward the chance floor as $\lambda_1$ increases, reaching
a median of $0.688$ over the members with $\lambda_1 \geq 5$. Broadly speaking, the model
reproduces $\mu_\nu$ more closely for the more chaotic dynamics. 
Where $\arf = 1$, the two samples are separable without error, and
Figure~\ref{fig:qualitative} shows what this looks like: at a member whose
trajectory is steady in time, the held-out states form a narrow bundle about
a single profile while the model produces states of the correct spatial scale
spread broadly in phase. 

The two panels differ in how orderly the trend is. Viewed against $L$, accuracy
tends to improve as the domain lengthens but is interrupted throughout by members
that revert to $\arf = 1$; these are the windows of non-chaotic behavior
that occur within predominantly chaotic ranges of the parameter, visible in
the right panel as dark points scattered along the top edge at large $L$.
Viewed against $\lambda_1$, the same points collapse into a single clear trend,
with the interruptions gathering at the low-$\lambda_1$ end.
This is the sense in which not the control parameter $\nu$ but $\lambda_1$ 
best describes the family. Accuracy first departs from $1$ at $\lambda_1 \approx
5$, which is what motivated the threshold defining the strongly chaotic
sub-family; the residual scatter above it narrows as $\lambda_1$ grows
further.

Restricting training to that sub-family improves the fit on its members.
$63$ cases survive the filter.
The filtered model attains a median $\arf$ of $0.614$ against the full model's
$0.688$ at the same test values, reducing the excess over chance by roughly
$40\%$. Training on the weakly chaotic members therefore costs accuracy on reproducing
the strongly chaotic ones, even though the two models are otherwise
identical. Throughout, the test values interleave the training grid rather
than lying beyond it (Appendix~\ref{app:solver}), so these results characterise interpolation within the
family and not extrapolation outside it.

\subsection{Spectral and physical consistency}
\label{sec:physics}

\begin{figure}[t]
  \centering
  \begin{subfigure}[t]{0.48\textwidth}
    \includegraphics[width=0.86\textwidth]{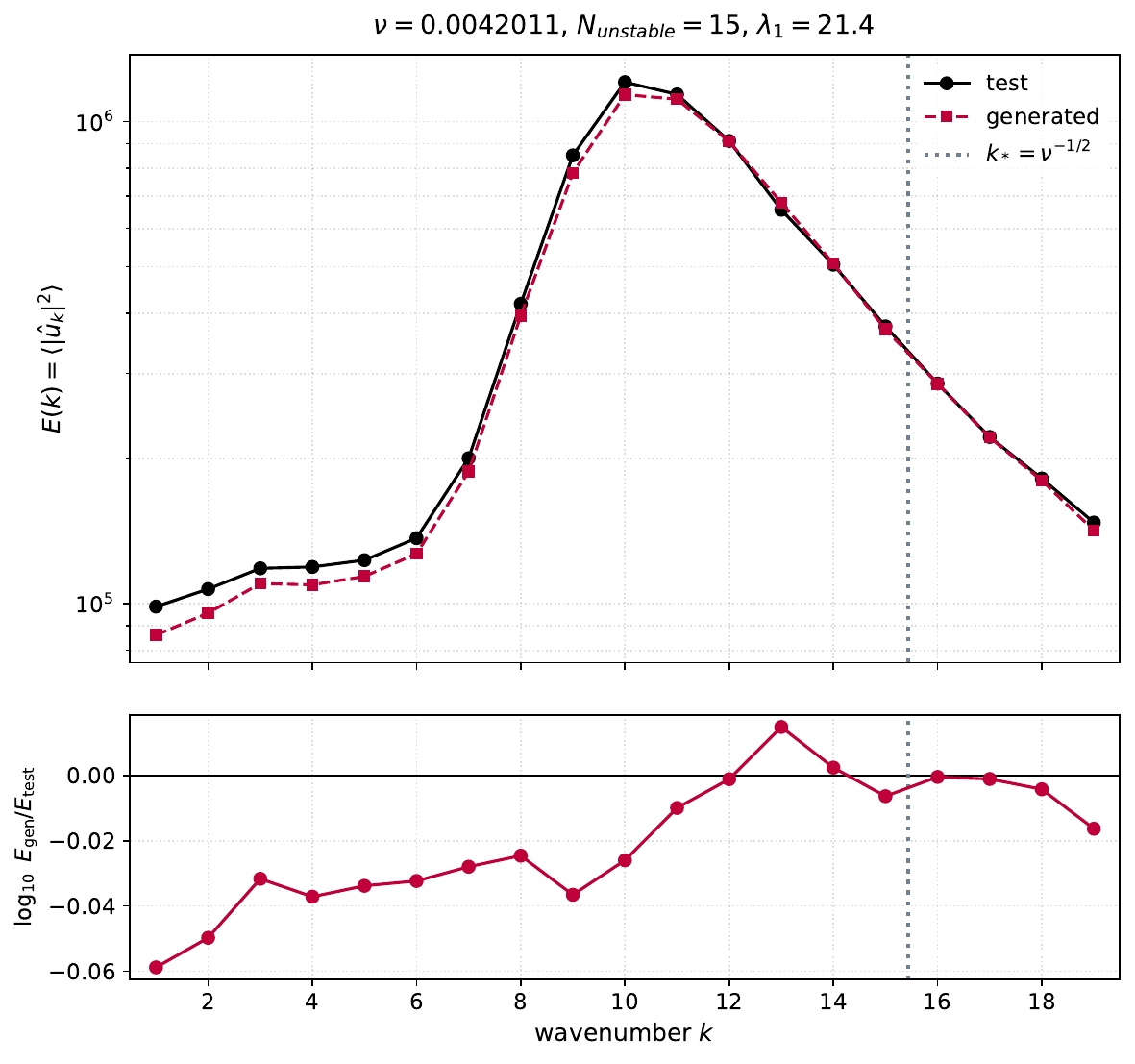}
    \caption{}
    \label{fig:spectrum-overlay}
  \end{subfigure}
  \hfill
  \begin{subfigure}[t]{0.48\textwidth}
    \includegraphics[width=\textwidth]{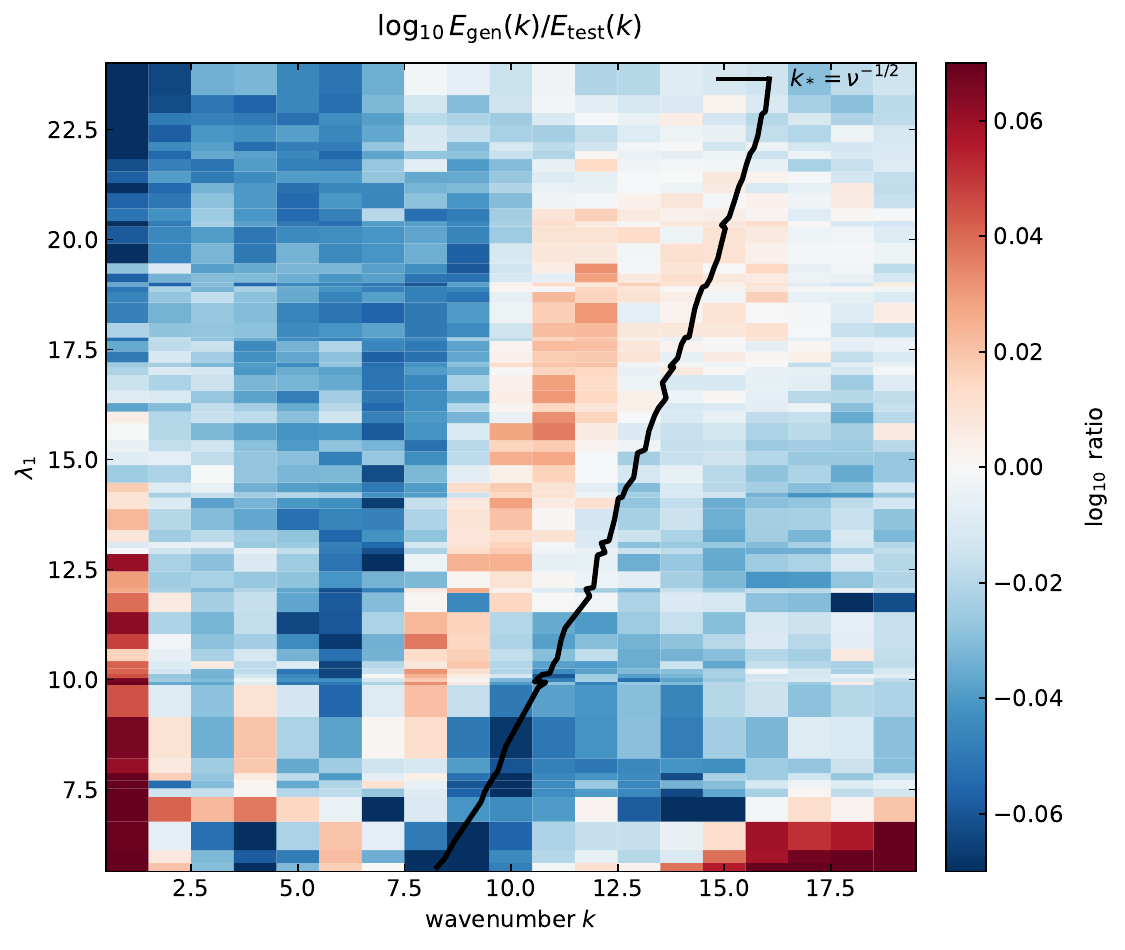}
    \caption{}
    \label{fig:spectrum-heatmap}
  \end{subfigure}
  \caption{Energy spectra of generated and held-out samples for
  the model trained on the $\lambda_1 \geq 5$ sub-family.
  \subref{fig:spectrum-overlay}: Mean modal energies
  $E_k = \langle |\hat{u}_k|^2 \rangle$ at a single strongly chaotic test
  value ($\nu = 0.00420$, $\lambda_1 = 21.4$), with the per-mode log ratio
  below. \subref{fig:spectrum-heatmap}: The same log ratio for every test
  value, against wavenumber and $\lambda_1$. The dotted line and the solid
  curve mark the neutral wavenumber $k_* = \nu^{-1/2}$.}
  \label{fig:spectrum}
\end{figure}

\begin{figure}[t]
  \centering
  \includegraphics[width=\textwidth]{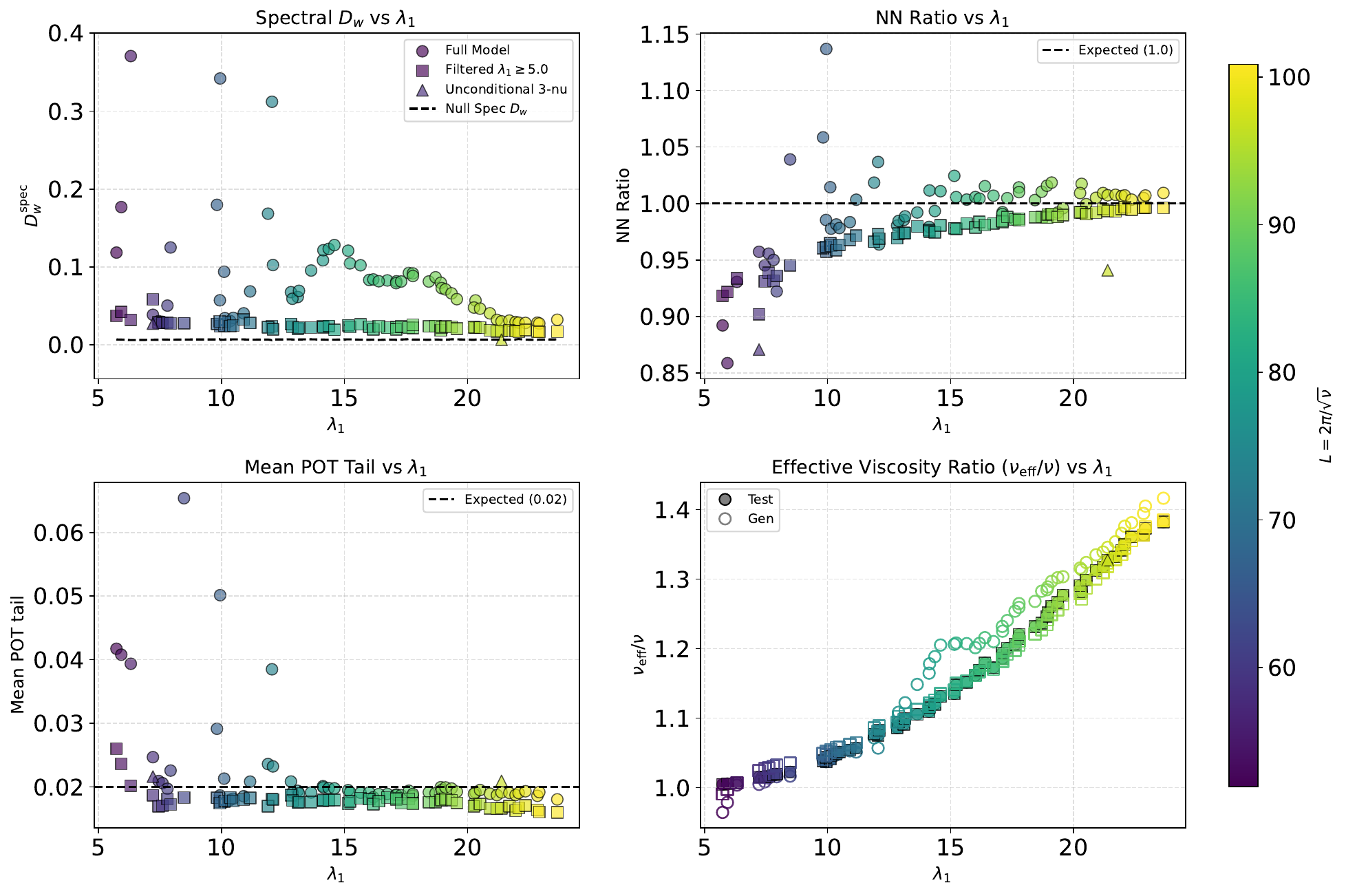}
  \caption{Supporting metrics of Section~\ref{sec:evaluation} against
  $\lambda_1$, for the members with $\lambda_1 \geq 5$; the weakly chaotic
  members are covered in Figure~\ref{fig:rfaccuracy}. Markers as in
  Figure~\ref{fig:rfaccuracy}. Dashed lines give the value attained by a
  model reproducing $\mu_\nu$ exactly: the resampling floor for $D_w$, unity
  for $\rho_{\mathrm{NN}}$, and $2\alpha = 0.02$ for tail coverage. In the
  lower right panel, filled markers are the held-out estimate of
  $\nu_{\mathrm{eff}}$ and open markers the generated estimate. One point of the full-family model, at $\lambda_1 \approx 8.5$ with $D_w \approx 1.0$, lies above the range of the upper-left panel and is not shown.}
  \label{fig:metrics}
\end{figure}

Figure~\ref{fig:spectrum} compares the generated and held-out energy spectra.
The generated spectrum reproduces the shape of the true one, including the
location of its peak and the decay above it, and the residual discrepancy is
systematic rather than diffusive: the log ratio is negative throughout the
energy-containing range below the peak, from about $-0.06$ at $k = 1$ to
$-0.03$ near $k = 8$, and crosses zero at the peak itself. The model
therefore holds too little energy at large scales, by roughly $8$--$13\%$ of
the modal energy, with the shortfall largest at the lowest wavenumbers. The
heatmap shows the same deficit across the family, together with a weak
surplus at wavenumbers just above the peak that follows the neutral
wavenumber $k_* = \nu^{-1/2}$ as the family is traversed; the redistribution
is tied to the scale at which the linear dynamics change character rather
than to a fixed wavenumber.

The supporting metrics agree with the classifier, and Figure~\ref{fig:metrics}
reports them over the strongly chaotic cases. For the filtered model, the
spectral discrepancy $D_w$ sits on the resampling floor across the whole
range, the nearest-neighbor ratio lies within a few percent of unity, and
tail coverage is within a few thousandths of $0.02$. Above
$\lambda_1 \approx 15$, two of these depart from their reference values in a
consistent direction: $\rho_{\mathrm{NN}}$ falls slightly below $1$ and tail
coverage slightly below $0.02$, indicating a small amount of generated mass
placed away from the attractor together with marginally thin tails. The
full-family model instead shows a pronounced hump in $D_w$ between
$\lambda_1 \approx 13$ and $18$, reaching about $0.10$, and the same members
carry the largest gap between the generated and held-out effective
hyperviscosities. One member near $\lambda_1 \approx 8.5$ is fit poorly by the
full-family model on both $D_w$ and $\arf$, and is fit normally by the
filtered model.

The effective hyperviscosity provides an independent check on the same samples.
Generated and held-out estimates of $\nu_{\mathrm{eff}}$ agree closely
throughout, with a median relative difference of $0.2\%$ and a worst case
below $2\%$, so the generated measure satisfies the energy balance to the
accuracy with which the held-out data does. Both estimates rise above $\nu$
as $\lambda_1$ increases, reaching some $40\%$ above it at the most chaotic
members. This is the truncation effect anticipated in
Section~\ref{sec:evaluation}: as $\nu$ decreases, the neutral wavenumber moves
toward the state mode cutoff $K$, and the discarded modes contribute more to
$\sum_k k^4 E_k$ than to $\sum_k k^2 E_k$, biasing the ratio upward. It
affects the two samples equally and does not enter the comparison between
them.
The size of this bias does not indicate that the state representation is too
coarse. The cutoff $K$ was chosen so that every linearly unstable mode is
retained; across the strongly chaotic members, the neutral
wavenumber $k_*$ runs from roughly $8.8$ to $16.1$, below $K = 19$
throughout. What the bias does reflect is that $\nu_{\mathrm{eff}}$ is by construction
the most truncation-sensitive quantity we compute, since its denominator
weights each mode by $k^4$ and is therefore dominated by the smallest
retained scales, which lie next to the modes the truncation discards.
The discarded modes are not energetically
negligible for that ratio, even though they are damped. The adequacy of the
representation for the learning problem is established instead by the
metrics that act on the retained modes, all of which agree that the generated
and held-out distributions coincide there.

\subsection{The cost of amortization}
\label{sec:amortization}

The conditional models are trained once and cover the whole family; the
alternative is to train a separate model for each parameter value. To measure
what the first choice costs, we trained three unconditional models, each on
the samples at a single $\nu$, with architecture and optimization identical to
the conditional models apart from the removal of the conditioning input. The
three values were selected before the runs from the test set, one from the
weakly chaotic regime and two from the strongly chaotic one. Their accuracies
appear as triangles in Figure~\ref{fig:rfaccuracy} and Figure 4, and samples from them appear in
Figure~\ref{fig:qualitative}. In the weakly chaotic case, the unconditional
model reaches $\arf = 0.997$, essentially the same complete separability as
the full-family model: the difficulty there is a property of the target
measure and not of amortization. In the two chaotic cases, the single-$\nu$
models do not improve on the conditional ones: at $\lambda_1 \approx 7.2$ the
unconditional model gives $\arf = 0.859$ against $0.775$ for the full-family
model and $0.720$ for the filtered one, and at $\lambda_1 \approx 21.4$ it
gives $0.639$ against $0.685$ and $0.600$ respectively.

Amortization therefore has no measurable impact on this family. At the
moderately chaotic member, it is better than free: a single model trained
across $160$ parameter values fits that member more closely than a model
trained on it alone, so pooling data across the family transfers rather than
interferes. The filtered conditional model is the most accurate of the three
at both chaotic members. 

\section{Discussion}
\label{sec:discussion}

We have shown that a single conditional diffusion model, trained once,
reproduces the invariant measure of the Kuramoto--Sivashinsky equation across
the strongly chaotic portion of a one-parameter family, with an accuracy that improves with the strength of the chaos over most of that range.
The model trained on that
sub-family is more accurate at both chaotic test cases than an unconditional
model trained on that case alone, so the amortization across $160$ parameter
values is obtained at no cost in fidelity. Our evaluation takes its starting
point from \citet{finn2024}: we retain several of their tests, add others
suited to the spatial structure of a PDE, and read the results against the
values an exact model would attain and against comparisons within the family.

The results also suggest that the feasibility of this approach is determined by the target rather than by the training. For a case
whose trajectory is steady in time or is a traveling wave, the invariant
measure is supported on the translation orbit of a single profile, a
one-dimensional curve in $\mathbb{R}^d$.
A generative diffusion model of the kind used here produces samples whose
distribution has full support in $\mathbb{R}^d$ and cannot concentrate on a
curve \citep{pidstrigach2022}; the separability reported in
Section~\ref{sec:headline} is the consequence.
The measure acquires the
full-dimensional support a diffusion model must produce only as the dynamics
become chaotic, so the argument identifies in advance, from the dynamics
alone, which cases of a parametrized family a model of this class can represent.

\section*{AI use statement}

In this work, we used generative AI tools (Claude Fable 5, Claude Opus 5 and
Gemini 3.6) for designing experiments, implementing methods, formulating
mathematical claims, interpreting results and qualitative analysis. We have
not used generative AI tools for generating synthetic datasets, which are
produced by a deterministic spectral solver, or for refining hypotheses.
Developing theoretical models, assisting with proofs, translation assistance
and data cleaning are not applicable to this work. Additionally, we used
generative AI tools for drafting the manuscript, analyzing the literature.

All algebra was verified independently before inclusion, and code was verified by
the two authors. We take responsibility for the final content of this work,
including text, claims or artifacts produced with the aid of generative AI.

\section*{Reproducibility statement}

The problem setup, including the parameter family and its division into
training, test and validation values, is given in Sections~\ref{sec:setup}
and~\ref{sec:data}, and the derivations underlying the evaluation are given in
Appendix~\ref{app:derivations}. Appendix~\ref{app:solver} specifies the
solver, its resolution, the time grid and the retained window of each
trajectory; Appendix~\ref{app:training} specifies the diffusion model, the
network and the training hyperparameters, together with the library version
used; and Appendix~\ref{app:metrics} specifies the classifier two-sample
test. The remaining metrics are defined in full in
Section~\ref{sec:evaluation}. All random draws in data generation, training,
sampling and evaluation are made from fixed seeds.
The code for data generation, training and evaluation will be released
publicly with the final version of the paper.

\bibliographystyle{iclr2027_conference}
\bibliography{ks.bib}

\begin{appendices}

\section{Derivations}
\label{app:derivations}

\subsection{Equivalence of the two conventions}
\label{app:rescaling}

The KS literature more commonly fixes the coefficients and varies the domain,
writing
\begin{equation}
\label{eq:ksL}
v_\tau + v v_y + v_{yy} + v_{yyyy} = 0, \qquad y \in [0, L),
\end{equation}
with periodic boundary conditions. Set $a = 2\pi/L$ and let
$v(y,\tau) = a\,u(ay,\,a^{2}\tau)$, so that $y \in [0,L)$ corresponds to
$x = ay \in [0,2\pi)$. Then $v_\tau = a^{3}u_t$, $v_y = a^{2}u_x$,
$v v_y = a^{3}u u_x$, $v_{yy} = a^{3}u_{xx}$ and $v_{yyyy} = a^{5}u_{xxxx}$,
so that \eqref{eq:ksL} becomes, after division by $a^{3}$,
\begin{equation}
u_t + u u_x + u_{xx} + a^{2} u_{xxxx} = 0 ,
\end{equation}
which is \eqref{eq:ks} with $\nu = a^{2} = (2\pi/L)^{2}$, equivalently
$L = 2\pi\nu^{-1/2}$, as stated in \eqref{eq:L}. Under the same change of
variables a Fourier mode $e^{i(2\pi n/L)y}$ becomes $e^{inx}$, which is why
the wavenumbers retained in \eqref{eq:fourier} are integers.

\subsection{Conservation of the spatial mean}
\label{app:mean}

Integrate \eqref{eq:ks} over $[0,2\pi)$. The advective term, $u u_x = \partial_x(u^{2}/2)$, and the two linear terms,
$u_x$ and $u_{xxx}$, integrates to zero
because by periodic boundary. Hence
\begin{equation}
\frac{d}{dt}\int_{0}^{2\pi} u(x,t)\, dx = 0 ,
\end{equation}
so $\hat{u}_0 = (2\pi)^{-1}\int_0^{2\pi} u\,dx$ is constant in time. Our
initial condition has zero mean, so $\hat{u}_0 \equiv 0$ throughout; being
constant along every trajectory, it carries no dynamical information and is
excluded from the state \eqref{eq:state}.

\subsection{Lyapunov exponents}
\label{app:lyapunov}

\paragraph{The variational equation.}
Let $u(x,t)$ be a solution of \eqref{eq:ks}, and let $u(x,t) + w(x,t)$ be a
second solution of the same equation, so that
\begin{equation}
(u+w)_t + (u+w)(u+w)_x + (u+w)_{xx} + \nu (u+w)_{xxxx} = 0 .
\end{equation}
Expanding the advective term as
$(u+w)(u+w)_x = u u_x + u w_x + w u_x + w w_x$ and subtracting
\eqref{eq:ks}, which $u$ satisfies by assumption, leaves
\begin{equation}
\label{eq:exact-perturbation}
w_t + u w_x + w u_x + w w_x + w_{xx} + \nu w_{xxxx} = 0 .
\end{equation}
Equation \eqref{eq:exact-perturbation} is exact and holds for a perturbation
of any size. linearization consists of discarding the single nonlinear term
$w w_x$, which is quadratic in the perturbation and is therefore negligible
beside the remaining terms when the perturbation is small. Writing $\delta$
for a perturbation governed by the resulting linear equation, and combining
$u \delta_x + \delta u_x = (u\delta)_x$ by the product rule, gives the
variational equation
\begin{equation}
\label{eq:variational}
\delta_t + (u\delta)_x + \delta_{xx} + \nu\, \delta_{xxxx} = 0 .
\end{equation}
This is a linear equation in $\delta$ whose coefficients depend on the
particular solution $u$ about which it was derived, so the same initial
perturbation evolves differently along different trajectories. It governs
infinitesimal perturbations only. A perturbation of finite size obeys
\eqref{eq:exact-perturbation}, in which the discarded term ceases to be
negligible once $w$ has grown.

\paragraph{The Lyapunov spectrum.}
Equation \eqref{eq:variational} is linear in $\delta$, so an initial
perturbation $\delta(\cdot,0)$ determines $\delta(\cdot,t)$ at all later
times, and multiplying $\delta(\cdot,0)$ by a constant multiplies
$\delta(\cdot,t)$ by the same constant. The asymptotic growth rate
\begin{equation}
\label{eq:lyap}
\lambda\bigl(\delta(\cdot,0)\bigr) = \lim_{t \to \infty} \frac{1}{t}
\ln \frac{\lVert \delta(\cdot,t) \rVert}{\lVert \delta(\cdot,0) \rVert},
\qquad \lVert \cdot \rVert = \lVert \cdot \rVert_{L^{2}([0,2\pi))},
\end{equation}
is therefore unchanged when $\delta(\cdot,0)$ is rescaled: it depends on the
direction of the initial perturbation and not on its magnitude. Different
directions may grow at different rates, and the values taken by
\eqref{eq:lyap} as the initial direction varies are the Lyapunov exponents
of the system, conventionally ordered
$\lambda_1 \geq \lambda_2 \geq \cdots$. Under the ergodicity assumed in
Section~\ref{sec:measure}, they are properties of the attractor rather than
of the particular trajectory along which \eqref{eq:variational} was
integrated.

Those directions whose growth rate is smaller than $\lambda_1$ are confined
to a proper subspace of the space of perturbations. A perturbation chosen
without reference to the dynamics therefore has a nonzero component along the
fastest-growing directions, and that component comes to dominate the others,
so that \eqref{eq:lyap} returns $\lambda_1$ for all but an exceptional set of
initial perturbations. This is why $\lambda_1$, alone among the exponents, is
readily accessible to measurement, and it is the only one we use.

\paragraph{$\lambda_1$ as a measure of chaos.}
If $\lambda_1 \leq 0$, no direction of perturbation grows and states that
begin close together remain close for all time; the trajectory approaches a
fixed point, a traveling wave or a periodic orbit. If $\lambda_1 \ge 0$, an
infinitesimal perturbation along a generic direction is amplified by a factor
$e^{\lambda_1 t}$, so that differences too small to measure are magnified
into differences of order the size of the attractor within a finite time.
This sensitive dependence on initial conditions is the defining property of
chaotic dynamics.

The exponent also measures how strong that sensitivity is. The time
$1/\lambda_1$ is the interval over which an uncertainty is amplified by a
factor $e$, so an uncertainty of relative size $\epsilon$ grows to order
unity after a time of roughly $\lambda_1^{-1}\ln(1/\epsilon)$. A member of
the family with a larger $\lambda_1$ therefore destroys information about its
initial state faster, and can be predicted over a shorter horizon, than a
member with a smaller one. This is the sense in which we call it more
strongly chaotic, and the sense in which $\lambda_1$ orders the family in
Section~\ref{sec:results}.

Both statements concern \eqref{eq:variational} and hold only while the
perturbation remains small enough for the linearization to be valid. Two
genuinely distinct solutions separate exponentially at first, but the
separation saturates once it is comparable with the extent of the attractor,
after which they are no more distant than two independent draws from
$\mu_\nu$.

\paragraph{Units and comparison with published values.}
An exponent is a rate, so its numerical value depends on the time variable,
and the two conventions of Appendix~\ref{app:rescaling} use different ones.
There $t = \nu\tau$, and the norm ratio in \eqref{eq:lyap} is unaffected by
the constant amplitude factor relating $u$ and $v$, so an exponent
$\lambda_1^{(L)}$ measured per unit $\tau$ for \eqref{eq:ksL} corresponds to
$\lambda_1 = \lambda_1^{(L)}/\nu$ per unit $t$ for \eqref{eq:ks}. Published
values are therefore directly comparable with ours: \citet{edson2019} report
$\lambda_1^{(L)} \approx 0.09$ at $L \approx 100$, and our largest member has
$L = 101.16$, that is $\nu \approx 0.00386$, so the conversion gives
$\lambda_1 \approx 23.3$ against the $23.7$ we measure there.

\subsection{The energy balance and the effective hyperviscosity}
\label{app:nueff}

Multiply \eqref{eq:ks} by $u$ and integrate over $[0,2\pi)$. The first term
is $\int u u_t\,dx = \tfrac{d}{dt}\tfrac{1}{2}\int u^{2}dx$, and the advective
term vanishes by periodicity, since $u \cdot u u_x = \partial_x(u^{3}/3)$.
The remaining two are integrated by parts, the boundary terms vanishing for
the same reason:
\begin{equation}
\int_0^{2\pi} u\, u_{xx}\, dx = -\int_0^{2\pi} u_x^{2}\, dx ,
\qquad
\int_0^{2\pi} u\, u_{xxxx}\, dx = \int_0^{2\pi} u_{xx}^{2}\, dx ,
\end{equation}
the second after two applications. Collecting the four terms,
\begin{equation}
\label{eq:energy-rate}
\frac{d}{dt}\, \frac{1}{2}\int_0^{2\pi} u^{2}\, dx
= \int_0^{2\pi} u_x^{2}\, dx - \nu \int_0^{2\pi} u_{xx}^{2}\, dx ,
\end{equation}
so the second-derivative term of \eqref{eq:ks} injects energy and the
fourth-derivative term removes it. By Parseval these two integrals equal
$2\pi\sum_k k^{2}\lvert\hat{u}_k\rvert^{2}$ and
$2\pi\sum_k k^{4}\lvert\hat{u}_k\rvert^{2}$, with the same constant in each.

Now take expectations with respect to $\mu_\nu$. Because $\mu_\nu$ is
invariant under the dynamics, the expectation of any observable is
independent of time, and the expectation of its time derivative therefore
vanishes. Applied to $\tfrac{1}{2}\int u^{2}dx$ through
\eqref{eq:energy-rate}, this gives
\begin{equation}
\label{eq:balance}
\sum_{k} k^{2} E_k = \nu \sum_{k} k^{4} E_k ,
\qquad E_k = \bigl\langle \lvert \hat{u}_k \rvert^{2} \bigr\rangle ,
\end{equation}
which rearranges to \eqref{eq:nueff}.
The sums run over all $k$; since $\hat{u}_{-k} = \overline{\hat{u}_k}$,
restricting them to $k \geq 1$ halves both sides and leaves their ratio
unchanged.

Read as a statement about the measure, \eqref{eq:balance} determines the
parameter of the equation from the modal energies alone:
\begin{equation}
\label{eq:nu-from-balance}
\nu = \frac{\sum_k k^{2} E_k}{\sum_k k^{4} E_k} ,
\end{equation}
provided the $E_k$ are those of the stationary measure and the sums are
complete. The effective hyperviscosity \eqref{eq:nueff} is the right-hand
side of \eqref{eq:nu-from-balance} evaluated on a finite sample and over the
retained modes $k \leq 19$ alone. It is called effective rather than simply
$\nu$ because neither hypothesis need hold: it returns the parameter of the
equation when the sample is drawn from $\mu_\nu$ and the sums are complete,
and otherwise returns whatever balance the sample does satisfy. This is what
makes it a diagnostic rather than a definition, condensing the modal energies
into a single number whose correct value is fixed by the dynamics instead of
being estimated from the data.

\section{Implementation details}
\label{app:implementation}

\subsection{Data generation}
\label{app:solver}

Trajectories are generated with a Fourier pseudospectral discretisation of
\eqref{eq:ks} retaining $100$ modes, advanced by the fourth-order exponential
time-differencing Runge--Kutta scheme ETDRK4 \citep{kassam2005}, which treats
the stiff linear part exactly. The solver resolution is therefore
considerably finer than the $K = 19$ modes retained in the state
\eqref{eq:state}.

Each trajectory is integrated to $T_{\mathrm{fin}} = 2000$ in $m = 110001$
steps of size $\tau = T_{\mathrm{fin}}/m \approx 0.018182$. The solver
discards the first $100001$ steps and stores the remaining $10000$ snapshots,
which are retained for training and evaluation. These are
global steps $100002$ through $110001$, spanning $181.8$ time units at the
end of the trajectory, and all of them are used without subsampling. The
discarded portion amounts to $1818.2$ time units.

All random draws in the generation, training, sampling, and evaluation pipelines are made
from fixed seeds.

\subsection{The diffusion model}
\label{app:training}

All models are variance-preserving diffusion models built with BayesFlow
2.0.12 \citep{radev2023} on a JAX backend. The forward process noises a state
$\uhat$ according to
\begin{equation}
\label{eq:forward}
\uhat_s = \alpha(s)\,\uhat + \sigma(s)\,\boldsymbol{\varepsilon},
\qquad \boldsymbol{\varepsilon} \sim \mathcal{N}(0, I_d),
\qquad s \in [0,1],
\end{equation}
with signal and noise scales set by the log signal-to-noise ratio
$\lambda(s)$ through $\alpha = \operatorname{sigmoid}(\lambda)^{1/2}$ and
$\sigma = \operatorname{sigmoid}(-\lambda)^{1/2}$. The cosine schedule \citep{nichol2021} sets
\begin{equation}
\lambda(s) = -2 \ln \tan\!\left(\frac{\pi \tilde{s}}{2}\right),
\qquad
\tilde{s} = \tilde{s}_{\min} + (\tilde{s}_{\max} - \tilde{s}_{\min})\, s ,
\end{equation}
the endpoints being fixed by the inverse relation
$\tilde{s} = (2/\pi)\arctan(e^{-\lambda/2})$ evaluated at $\lambda = 15$ and
$\lambda = -15$, which truncates the log signal-to-noise ratio to
$[-15, 15]$.

The network predicts the velocity
$\mathbf{v}_s = \alpha(s)\boldsymbol{\varepsilon} - \sigma(s)\uhat$ \citep{salimans2022}, and the
objective is the squared error in that prediction, with $s$ drawn uniformly
on $[0,1]$ and each term weighted by
$w(\lambda) = \operatorname{sigmoid}(-\lambda + 2)$
\citep{kingma2023}; this is the
weighting the schedule supplies by default, and it is defined for the
noise-prediction form of the objective rather than the velocity form used
here. Sampling integrates the reverse-time stochastic differential equation
from $s = 1$ to $s = 0$ with a two-step adaptive stochastic integrator taking
an adaptive number of steps. No conditioning dropout is applied.

The velocity network is the time-conditioned residual multilayer perceptron
supplied by BayesFlow, with hidden widths $(256, 512, 1024, 1024, 512, 256)$,
mish activations, He-normal initialization, layer normalization, dropout
$0.05$ and a residual connection around each block; spectral normalization is
not used. The noised state and the conditioning scalar are each projected
linearly and merged by concatenation at the input. The diffusion time enters
separately, embedded in $128$ random Fourier features with scale $30$ and
injected at every hidden block by feature-wise linear modulation
\citep{perez2018}, which here learns an additive shift per feature and no
multiplicative scale.

States and the conditioning scalar are standardized using statistics of the
training set, and generated samples are mapped back to physical scale before
any evaluation.

Training uses Adam \citep{kingma2015} with a learning rate of
$10^{-3}$, batch size $512$, for $10$ epochs. The $16$ validation values of
Section~\ref{sec:data} were used to monitor the validation loss during
training; no early stopping or checkpoint selection was performed on them, so
the models reported are those obtained at the end of the tenth epoch.

\subsection{Classifier settings}
\label{app:metrics}

The two-sample classifier of Section~\ref{sec:evaluation} is a random forest
of $300$ trees with a minimum of five samples per leaf, as implemented in
scikit-learn. For each test value of $\nu$ the generated and held-out samples
are first balanced to a common size, which with $10{,}000$ states in each
retains all of them and makes $\arf = 0.5$ the exact chance baseline. The
pooled $20{,}000$ states, labelled by their origin, are split at random into
$70\%$ for fitting and $30\%$ for evaluation, stratified by label, and $\arf$
is the accuracy on the $6{,}000$ evaluation states.

No standardization is applied before fitting. The splits of a random forest
are thresholds on individual coordinates and are therefore invariant to the
scale of each coordinate, which matters here because the components of
\eqref{eq:state} span several orders of magnitude in variance across
wavenumbers. The constraint on leaf size is the only regularization used; it
limits how finely the forest can partition the training split, which is what
keeps the held-out accuracy an estimate of distinguishability rather than of
memorization.
\end{appendices}
\end{document}